# Fast Intrinsic Mode Decomposition of Time Series Data with Sawtooth Transform

Louis Yu Lu, CGBU, Oracle Corporation
E-mail: yu.lu@oracle.com

Abstract: An efficient method is introduced in this paper to find the intrinsic mode function (IMF) components of time series data. This method is faster and more predictable than the Empirical Mode Decomposition (EMD) method devised by the author of Hilbert Huang Transform (HHT). The approach is to transforms the original data function into a piecewise linear sawtooth function (or triangle wave function), then directly constructs the upper envelope by connecting the maxima and construct lower envelope by connecting minima with straight line segments in the sawtooth space, the IMF is calculated as the difference between the sawtooth function and the mean of the upper and lower envelopes. The results found in the sawtooth space are reversely transformed into the original data space as the required IMF and envelopes mean. This decomposition method process the data in one pass to obtain a unique IMF component without the time consuming repetitive sifting process of EMD method. An alternative decomposition method with sawtooth function expansion is also presented.

## 1    HHT Introduction

The real time series numeric data from natural phenomena, life science, social and economic system are mostly nonlinear and non-stationary. The traditional analysis methods like Fourier transform and Wavelet transform presume the linear and stationary of the underlying system, the predefined function bases have little bearing on the physical meaning of the system, and are difficult to reveal the nature of the real data. The recent researches tend to use adaptive function bases, Huang et al introduced the IMF as adaptive a posteriori function base in the form of Hilbert spectrum expansion, which has meaningful instantaneous frequency[1][2]. The HHT has broad application to signal processing and times series data analysis in the fields of health, environmental, financial and manufacturing industries [4].

### 1.1   Hilbert Transform

For an arbitrary time series $X(t)$, its Hilbert transform $Y(t)$ is defined as:

$$Y(t) = \frac{1}{\pi} P \int_{-\infty}^{\infty} \frac{X(\tau)}{t - \tau} d\tau \tag{1.1}$$

Where $P$ indicates the Cauchy principal value. With this definition, $X(t)$ and $Y(t)$ form the complex conjugate pair, so an analytic signal $Z(t)$ is defined as:

$$Z(t) = X(t) + iY(t) = a(t)e^{i\theta(t)} \tag{1.2}$$

In which:

$$a(t) = \sqrt{X^2(t) + Y^2(t)}, \quad \theta(t) = \arctan\frac{Y(t)}{X(t)} \qquad (1.3)$$

Theoretically, there are infinitely many ways of defining the imaginary part, but the Hilbert transform provides a unique way of defining the imaginary part so that the result is an analytic function. A comprehensive coverage of Hilbert transformation can be found in reference [3].

The instantaneous frequency is defined as:

$$\omega(t) = \frac{d\theta(t)}{dt} \qquad (1.4)$$

With this definition, not every function *X(t)* has well defined physical meaning [1]. Only a special family of function has the meaningful instantaneous frequency.

## 1.2   Intrinsic Mode Function (IMF)

Dr Huang analyzes the requirement of meaningful instantaneous frequency on Hilbert transformation, and introduced the IMF [1] [2], which satisfies the two conditions:
  1)  In the whole data set, the number of zero crossings must equal or differ at most by one;
  2)  At any point, the mean value of the envelope defined by the local maxima and the envelope defined by the local minima is zero.

Apparently the first condition is necessary for oscillation data; the second condition requires the symmetric upper and lower envelopes of IMF, as the IMF component is decomposed from the original data; it is quite challenging to find the real envelopes because of nonlinear and non-stationary nature in the data. Only a few functions have the known envelopes, for example, the constant amplitude sinusoidal function.

All the IMF has the meaningful instantaneous frequency defined by (1.4).

## 1.3   Empirical Mode decomposition (EMD)

To break up the original data into a series of IMF, Huang et al invented the EMD method [1] [4], the idea is to separate the data into a slow varying local mean part and fast varying symmetric oscillation part, the oscillation part becomes the IMF and local mean the residue, the residue serves as the input data again for further decomposition, the process repeats until no more oscillation can be separated from the last residue.

On each step of the decomposition, because the upper and lower envelope is unknown initially, a repetitive sifting process is applied to approximate the envelopes and obtain the IMF and residue as the following:

1) Find the extrema (maxima and minima) of the input data $m_0(t) = f(t)$, and connect the maxima with Spline function to form the upper envelope $U_0(t)$, and connect the minima with Spline function to form the lower envelope $L_0(t)$. The fist envelopes mean is calculated as:

$$m_1(t) = \frac{U_0(t) + L_0(t)}{2}$$

And the first difference between the data and the mean:

$$h_1(t) = m_0(t) - m_1(t)$$

2) Take the previous difference $h_{i-1}(t)$ as input again, find the envelopes the same way as in step 1), and obtain the mean and the difference:

$$m_i(t) = \frac{U_{i-1}(t) + L_{i-1}(t)}{2}$$

$$h_i(t) = m_{i-1}(t) - m_i(t)$$

3) Repeat the step 2 until the stopping condition is satisfied. The condition may be the small difference between successive $h_i(t)$ and $h_{i-1}(t)$ or small value of envelope mean $m_i(t)$. The first IMF and the residue:

$$c_1(t) = h_n(t)$$

$$r_1(t) = f(t) - c_1(t)$$

4) Take residue $r_1(t)$ as the input data again and repeat the process starting at step 1) until the number of extrema is less than 2.

The original data can be reconstructed by summing up all the IMFs and the last residue:

$$f(t) = \sum_1^n c_i(t) + r_n(t)$$

There are three problems with this decomposition method:
1) The Spline (cubic) connecting extrema is not the real envelope, some unwanted overshoot may introduced by the Spline interpolation, the resulting IMF function does not strictly guarantee the symmetric envelopes [1]. Higher order Spline does not in theory resolve the problem.
2) The different stopping condition value results different set of IMF [1], picking up the right value is a more subjective matter, making the results unpredictable.
3) The repetitive sifting process is time consuming, because the problem item 2), more sifting does not produce better results. Some researches try to improve the performance, the study in reference [6] only speedups the envelopes building process; it still requires the repetitive sifting process of EMD.

## 2    Sawtooth transform method

To overcome difficulty of finding the upper and lower envelopes in the original data space, alternative spaces are explored where the envelopes are easy to find, and the points between the original and the alternative spaces must have one to one mapping. The sawtooth space described below meets those requirements.

### 2.1    Transform time series data to sawtooth function

The required sawtooth (triangle wave) function can be constructed by connecting the successive extrema of the original data function with straight line segment, those extrema are alternating maxima and minima on the original data.

On each segment, the two extremities coincides the maxima and minima of the original data, in between, the variation of the original data value is monotone, it has one on one mapping to the straight line segment of the sawtooth function.

The original time series data *f(t)* has **m** extrema:

$$E(t_j) \qquad t_0 \le t_j \le t_{m-1}$$

There are *k* maxima among **m** extrema on the upper envelope *U(t)*:

$$E(t_{j-1}) < U(t_i) = E(t_j) > E(t_{j+1}) \quad 0 \le i < k$$

And *l* minima among **m** extrema on the lower envelope *L(t)*:

$$E(t_{j-1}) > L(t_i) = E(t_j) < E(t_{j+1}) \quad 0 \le i < l$$

Where *k* is either equals to *l* or differs by one, and

$$k + l = m$$

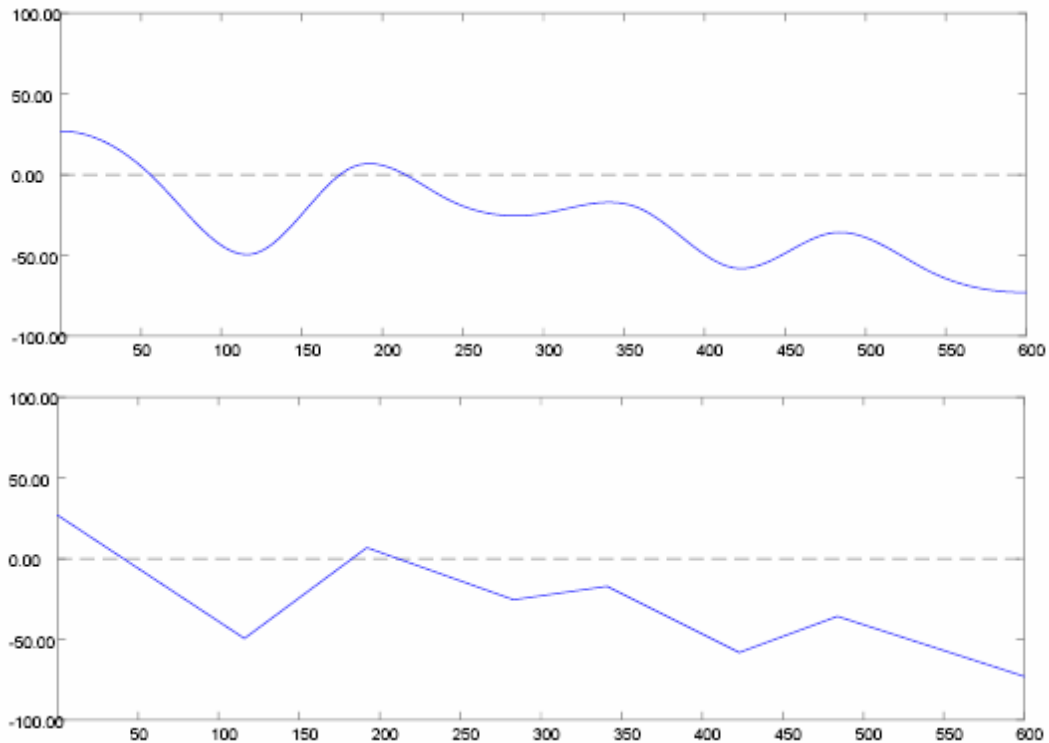

Figure 2.1 Original data and corresponding sawtooth (triangle wave) function

The sawtooth function is defined with **m-1** segments as:

$$s(t) = E(t_i) + (E(t_{i+1}) - E(t_i))\frac{t - t_i}{t_{i+1} - t_i} \quad 0 \leq i < m-1 \qquad (2.1)$$

Besides the points on the original data function, all the points in the time interval is applied the same transform. In the data space, a point is designated by **(t, x)**, while in the sawtooth space, a point has the coordinate **(u, s)**. The sawtooth transform on each segment is defined as:

$$u(t) = t_i + \frac{f(t) - E(t_i)}{E(t_{i+1}) - E(t_i)}(t - t_i) \quad t_i \leq t \leq t_{i+1} \quad 0 \leq i < m-1 \quad (2.2)$$

$$s(u) = x(t) \qquad (2.3)$$

The effect of this transform does not change the value on the vertical direction (signal value), it compress or expand the space on the horizontal direction (time) to make the data function into piecewise linear sawtooth function.

## 2.2 Finding envelopes, mean and IMF in sawtooth space

The sawtooth function varies linearly between alternating maxima **U(t_i)** and minima **L(t_i)**. The upper envelope can be constructed by connecting the successive maxima:

$$U(u) = U(u_i) + \frac{U(u_{i+1}) - U(u_i)}{u_{i+1} - u_i}(u - u_i) \quad 0 \leq i < k-1 \qquad (2.4)$$

And the similarly to build the lower envelope:

$$L(u) = L(u_i) + \frac{L(u_{i+1}) - L(u_i)}{u_{i+1} - u_i}(u - u_i) \quad 0 \leq i < l-1 \qquad (2.5)$$

The sum of maxima and minima equals the number of extrema:
The residue is the mean of the envelopes:

$$r(u) = \frac{U(u) + L(u)}{2} \qquad (2.6)$$

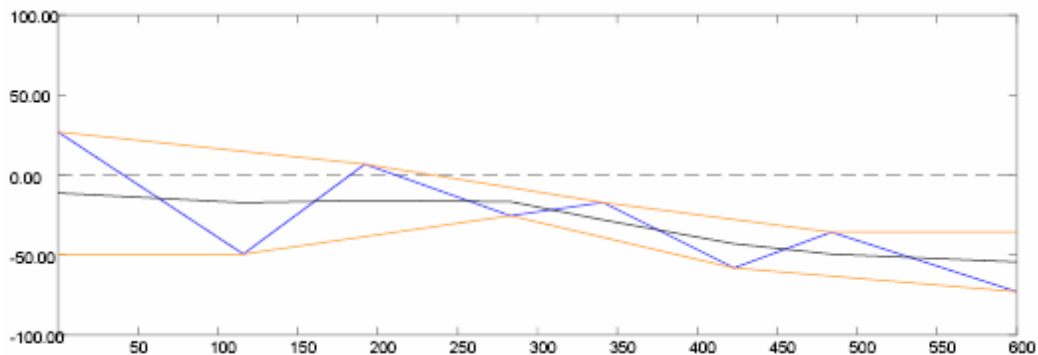

Figure 2.2 Upper and lower envelopes, envelopes mean and sawtooth function

The IMF is the difference between the sawtooth function and the mean:

$$c(u) = s(u) - r(u) \qquad (2.7)$$

The function **$c(u)$** has symmetric upper and lower envelope, it satisfies the two conditions of IMF specified in 1.2.

By subtracting the value of the mean from the upper and lower envelopes, they become symmetric around the value zero.

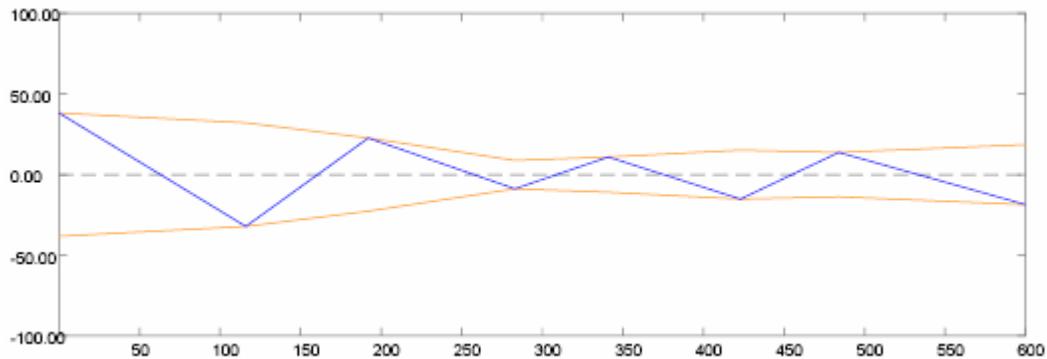

Figure 2.3 IMF and envelopes in the sawtooth space

In the sawtooth space, all the important functions: IMF, upper and lower envelopes are piecewise linear, the calculation is only required on the extrema. No repetitive process is required.

## 2.3   Transform back to original data space

To obtain the IMF **$c_{data}(t)$**, residue **$r_{data}(t)$** and envelopes in the original data space, it only needs to transform the results in the sawtooth space back to the original data space. With sawtooth transform, the value on vertical direction does not change, so the only requirement is to find the corresponding horizontal coordinate in the sawtooth space, which is defined by the transform (2.2). And those important reverse conversion functions are:

$$c_{data}(t) = c(u(t)) \qquad (2.9)$$
$$r_{data}(t) = r(u(t)) \qquad (2.10)$$
$$U_{data}(t) = U(u(t)) \qquad (2.11)$$
$$L_{data}(t) = L(u(t)) \qquad (2.12)$$

The reverse transform does not break the envelope symmetry of $c(u)$, so $c_{data}(t)$ still satisfies the two conditions of IMF.

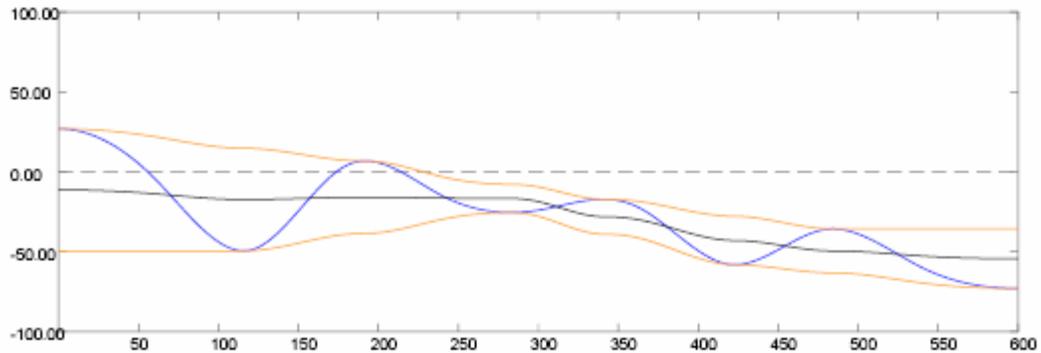

Figure 2.4 Upper and lower envelopes, envelopes mean and data in the original space

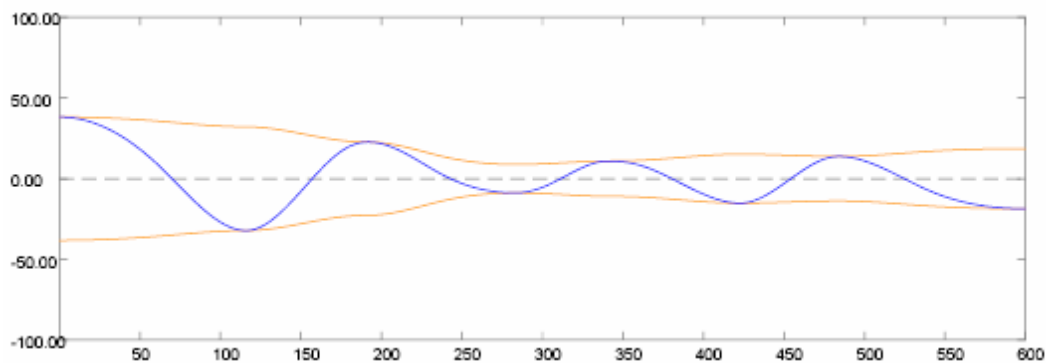

Figure 2.5 IMF and envelopes in the original space

## 2.4   IMF decomposition algorithm

With sawtooth transform, the intrinsic mode can be decomposed in the following steps

1) Find the maxima and minima of the input data. If the number of extrema is less than two, all the intrinsic modes are all found, stop the process.

2) Connect the extrema on the data with straight line segments to build the sawtooth function;

3) Connect the maxima with straight line segments to build the upper envelope in the sawtooth space, and connect the minima with straight line segments to build the lower envelope;

4) Use the formula (2.2) to calculate the $u$ coordinate in the sawtooth space for each time value $t$, calculate the residue with formula (2.6) and IMF with formula (2.7).

5) Transform the results back to the original space using formula (2.9), (2.10), (2.11) and (2.12).

6) Take the residue (envelopes mean) as input data again and back to step 1.

## 2.5   Function extension and boundary points

To calculate the envelopes and IMF on the boundary points, the data series must be extended on both ends to ensure the continuity on starting and ending points. Based on the nature of the data, there four possible extensions can be applied around the boundary points:

even extension; odd extension; cyclic extension can be used if the starting point and the ending point have the equal value; and the trend extension which extends the upper and lower envelopes on the first and last cycle.

The extension is only required to extend extrema in order to build the upper and lower envelopes in the sawtooth space; only two extra points are needed on each end.

1) Even extension: the starting point is considered as the first extrema $E(t_0)$, and the ending point as the last extrema $E(t_{m-1})$.

$$t_{-1} = t_0 - (t_1 - t_0) \qquad\qquad E(t_{-1}) = E(t_1) \qquad\qquad\qquad (2.13)$$
$$t_{-2} = t_0 - (t_2 - t_0) \qquad\qquad E(t_{-2}) = E(t_2) \qquad\qquad\qquad (2.14)$$
$$t_m = t_{m-1} + (t_{m-1} - t_{m-2}) \qquad E(t_m) = E(t_{m-2}) \qquad\qquad\quad (2.15)$$
$$t_{m+1} = t_{m-1} + (t_{m-1} - t_{m-3}) \qquad E(t_{m+1}) = E(t_{m-3}) \qquad\qquad\quad (2.16)$$

2) Odd extension: given the starting point $f(t_s)$, and the ending point $f(t_s)$.

$$t_{-1} = t_s - (t_0 - t_s) \qquad\qquad E(t_{-1}) = f(t_s) - (E(t_0) - f(t_s)) \qquad (2.17)$$
$$t_{-2} = t_s - (t_1 - t_s) \qquad\qquad E(t_{-2}) = f(t_s) - (E(t_1) - f(t_s)) \qquad (2.18)$$
$$t_m = t_e + (t_e - t_{m-1}) \qquad\quad E(t_m) = f(t_e) - (E(t_{m-1}) - f(t_e)) \qquad (2.19)$$
$$t_{m+1} = t_e + (t_e - t_{m-2}) \qquad\quad E(t_{m+1}) = f(t_e) - (E(t_{m-2}) - f(t_e)) \quad (2.20)$$

3) Cyclic extension: the starting point and end point are considered as extrema and they have the same value $E(t_0) = E(t_{m-1})$.

$$t_{-1} = t_0 - (t_{m-1} - t_{m-2}) \qquad E(t_{-1}) = E(t_{m-2}) \qquad\qquad\quad (2.21)$$
$$t_{-2} = t_0 - (t_{m-1} - t_{m-3}) \qquad E(t_{-2}) = E(t_{m-3}) \qquad\qquad\quad (2.22)$$
$$t_m = t_{m-1} + (t_1 - t_0) \qquad\quad E(t_m) = E(t_1) \qquad\qquad\qquad (2.23)$$
$$t_{m+1} = t_{m-1} + (t_{m-1} - t_{m-3}) \quad E(t_{m+1}) = E(t_2) \qquad\qquad\qquad (2.24)$$

4) Trend extension: by extending the first upper and lower envelopes to the left, the envelopes mean can be calculated on the starting point. Similarly, the envelopes mean on the ending point can be calculated by extending the last upper and lower envelopes to the right.

## 2.6  Flat maxima ceiling and minima floor

In some special case, the maxima is not a single point, it is instead a flat ceiling of constant value segment, in this case, the whole segment will be considered as the maxima when building the sawtooth function, both the left end and the right end of the segment connect to the descending contiguous segments. The maxima segment becomes a part of the upper envelope.

Similar for the minima floor consisting of constant value segment, both the left end and the right end of the segment connect to ascending contiguous segments, and the minima floor becomes part of the lower envelope.

# 3  Numerical example

The results in the following example are generated with a program coded in D language [5] for performance; the resulting data are converted into SVG format for display.

## 3.1 Original data, intrinsic modes and last residue

The simulated data is randomly generated. From the figures below, the last residue illustrates the general trend of the data. The IMF displays the fluctuation with different variation rate on different intrinsic mode.

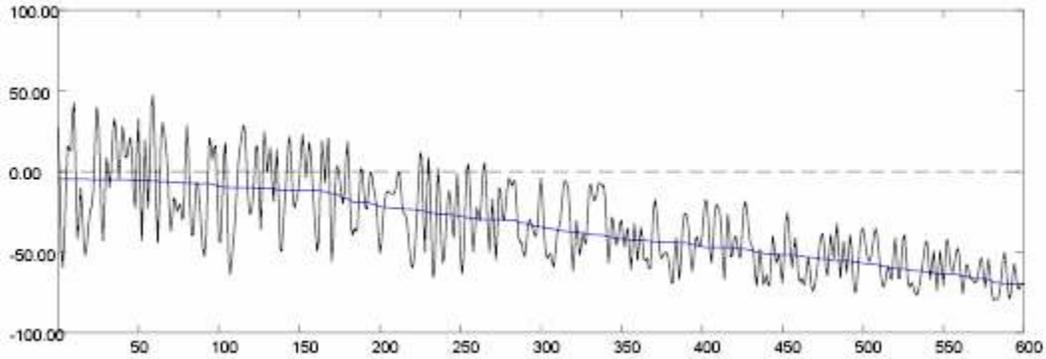

Figure 3.1 Original data and last residue

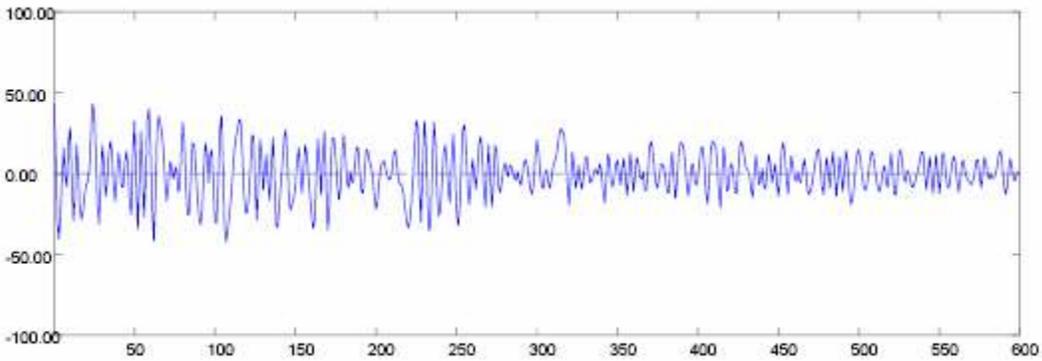

Figure 3.2 First IMF, 210 extrema

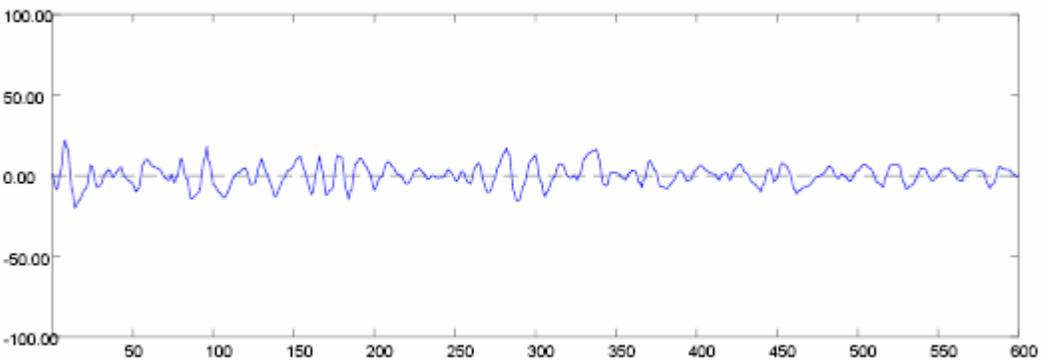

Figure 3.3 Second IMF, 85 extrema

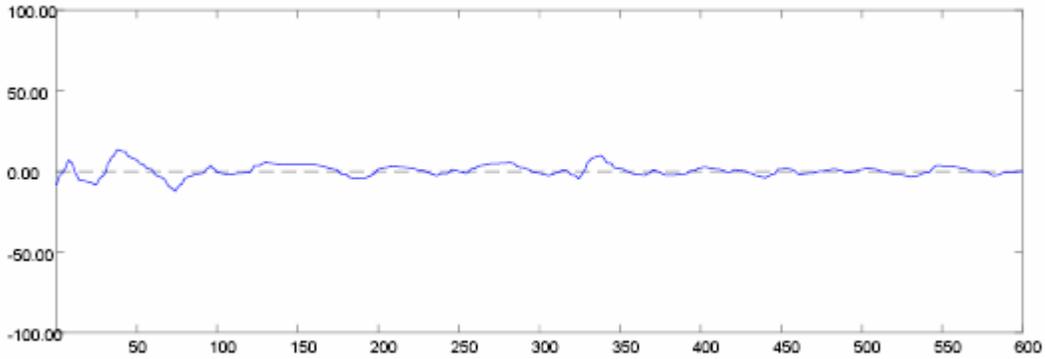

Figure 3.4 Third IMF, 33 extrema

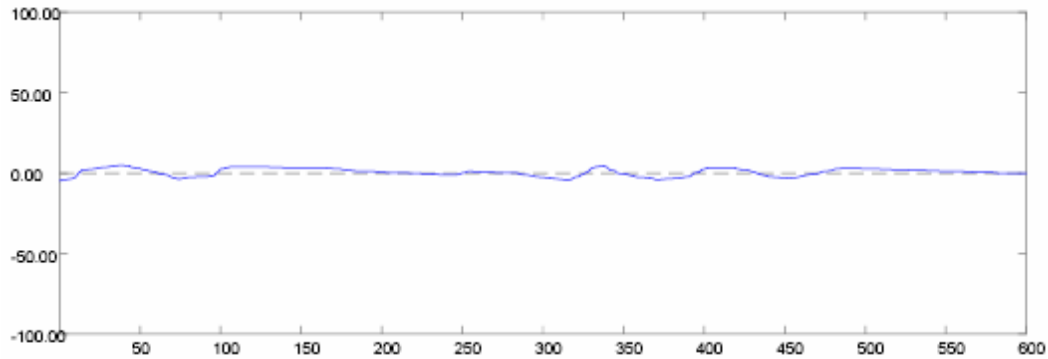

Figure 3.5 Fourth IMF, 12 extrema

## 3.2   Original data and envelopes means

The generated figures below display the residue (envelopes means) on different mode. The trend of the changes in different time scale is illustrated by the residue of different mode.

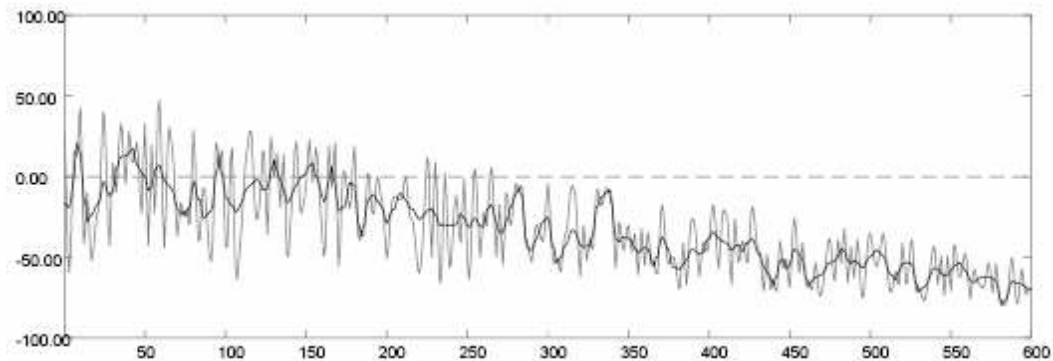

Figure 3.6 Data and first mode envelopes mean

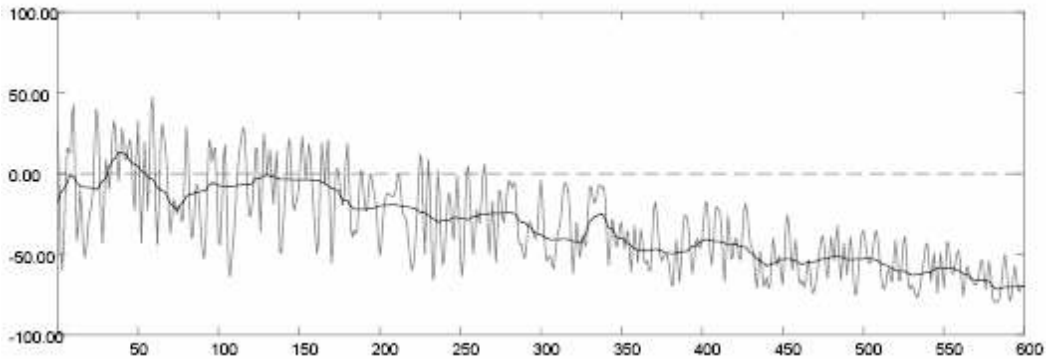

Figure 3.7 Data and second mode envelopes mean

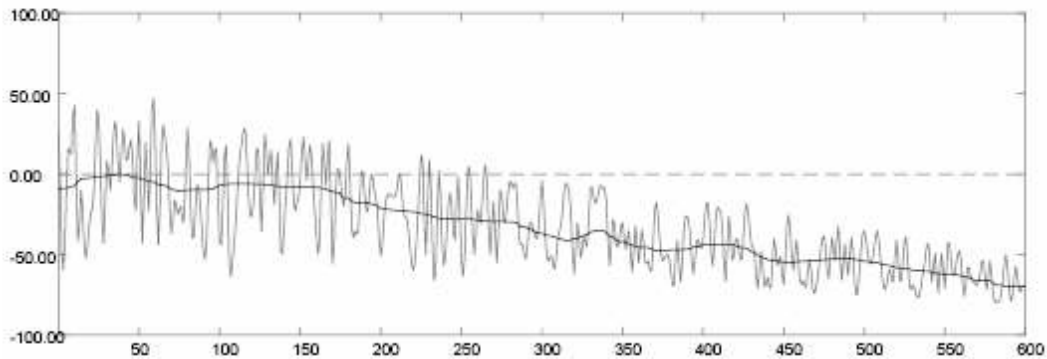

Figure 3.8 Data and third mode envelopes mean

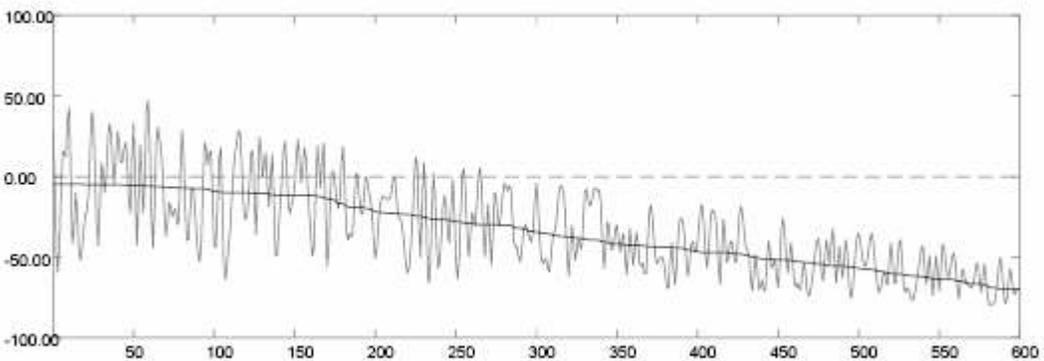

Figure 3.9 Data and fourth mode envelopes mean

## 4   IMF decomposition with sawtooth functions expansion method

An alternative approach can also be devised for IMF decomposition. In this method, the data is first expanded into a series of sawtooth functions with the following steps:

1) Take the original data as input;
2) Find the corresponding sawtooth function of the input data with the method described in section 2.1, and keep the sawtooth as one component;
3) Calculate the difference between the input data and the sawtooth function, and find the biggest value of the difference;
4) If the biggest difference value is less than a prescribed value, stop the expansion;
5) Take the difference as the input again, and go to step 2 for further processing.

With the sawtooth function components, the biggest difference between the original data and reconstructed data by summing up all the sawtooth functions will be less than the given value.

For each sawtooth function, its piecewise linear envelopes mean and IMF can be calculated with the method described in 2.2.

The envelopes mean of the original data will be the summation of all the envelopes means of each sawtooth component function with a maximal error less than the given value.

Similarly, The IMF of the original data will be the summation of all the IMFs of each sawtooth component function with a maximal error less than the given value.

This approach does not require the reverse transform in 2.3. It is less efficient compare with the method in section 2 for one dimensional data, but it may be easier to be generalized for two or higher dimensional data because only the extrema are considered for the processing, and each component is linear piecewise.

## 5   Discussion

A given time series data may have infinite number of IMF decompositions, but not all decomposition have good physical meaning [1].

Besides the envelopes mean method presented in this paper, other methods to calculate the residue function have been tested also, one is to connected the middle points between the neighboring sawtooth segment; the other one is connecting the weight center of the neighboring triangles formed by the consecutive three point on the sawtooth function, with those methods, the IMF and residue function are less smooth compare with the ones generated with envelopes mean method.

The algorithm of the sawtooth transform method only involves the data points spanned by the extent of four consecutive extrema which are the overlapping area of two maxima and two minima, only a small fraction of whole input data set is required for the running process. For streaming input data, high performance program can be implemented with small memory footprint. This feature can apply to the scenario where the data is collected in real time from a sensor; the processing can be started without waiting for the whole data set being collected.

## 6   Conclusion

The sawtooth transform method presented in this paper provides a fast and reliable way to decompose the time series data into a unique set of intrinsic mode components. It resolves the three problems related with the EMD method devised by the author of HHT.

This method will serve as a powerful tool to separate the nonlinear and non-stationary time series data into the trend part (envelopes mean) and the oscillation part (IMF) on different time scales, and will find application in many fields where traditionally Fourier analysis method or Wavelet analysis method dominate.